\documentclass{amsart}

\newtheorem{prop}{Proposition}[section]
\newtheorem{lemma}[prop]{Lemma}
\theoremstyle{definition}
\newtheorem{exa}[prop]{Example}

\usepackage{tikz-cd}
\usepackage{amssymb}

\newcommand{\grSpec}{\text{gr-Spec}}
\newcommand{\Spec}{\text{Spec}}

\title[Graded Injective Domains and Integral Extensions]{Graded Injective Domains, Valuation Overrings and Integral Extensions}

\author[Jha]{Ashish Jha}
\address{Indiana University East}
\email{jhaash@iu.edu}

\subjclass{13A02, 13B21, 13F05}
\keywords{graded integral domain, graded $i$-domain, injective extension, graded Pr\"ufer domain, graded valuation domain}

\begin{document}

\begin{abstract}
An extension $R \subseteq S$ of integral domains is an $i$-extension if the contraction map $\Spec(S) \to \Spec(R)$ is injective, and an integral domain $R$ is an $i$-domain if $R \subseteq S$ is an $i$-extension for every overring $S$ of $R$. Ira J. Papick proved the following two results, where $\overline{R}$ denotes the integral closure of $R$: $R$ is an $i$-domain if and only if $R \subseteq V$ is an $i$-extension for every valuation overring $V$ of $R$ and $R \subseteq \overline{R}$ is an $i$-extension; and if $T$ is an $i$-domain integral over $R$, then $R$ is an $i$-domain whenever $T$ is local or $R$ is integrally closed. In this paper, we prove graded analogues of these results. 
\end{abstract}

\maketitle

\section{Introduction}

Let $R$ be an integral domain with quotient field $K$. An \textit{overring} of $R$ is a ring $S$ with $R \subseteq S \subseteq K$. For a ring $S$, let $\Spec(S)$ denote the set of prime ideals of $S$. Papick \cite{Pap76} calls an extension of rings $R \subseteq S$ an \textit{injective extension} (or \textit{$i$-extension}) if the contraction map $\Spec(S) \to \Spec(R)$, given by $Q \mapsto Q \cap R$, is injective. Equivalently, for each prime ideal $P$ of $R$, at most one prime ideal $Q$ of $S$ satisfies $Q \cap R = P$. He calls $R$ an \textit{$i$-domain} if $R \subseteq S$ is an $i$-extension for every overring $S$ of $R$. Section 2 gives precise definitions of all the terms we use in this paper.

The \textit{integral closure} $\overline{R}$ of $R$ is the overring of $R$ consisting of the elements of $K$ that are roots of monic polynomials in $R[X]$. Recall that $R$ is a \textit{valuation domain} if $x \in R$ or $x^{-1} \in R$ for each nonzero $x \in K$. Papick shows that $R$ is a local $i$-domain if and only if $\overline{R}$ is a valuation domain \cite[Corollary~2.15]{Pap76}, and that $R$ is an $i$-domain if and only if $R \subseteq V$ is an $i$-extension for each valuation overring $V$ of $R$ and $R \subseteq \overline{R}$ is an $i$-extension \cite[Proposition~2.18]{Pap76}. He also shows that if $T$ is an $i$-domain that is integral over $R$, then $R$ is an $i$-domain whenever $T$ is local or $R$ is integrally closed \cite[Proposition~2.20]{Pap76}.

Many of the ideas described in \cite{Pap76} have graded analogues. A standard example of a graded domain is the polynomial ring $k[X_1, \ldots, X_n]$ over a field $k$, graded by total degree. Here, the \textit{total degree} of a monomial $X_1^{a_1} \cdots X_n^{a_n}$ is $a_1 + \cdots + a_n$, and a nonzero polynomial is \textit{homogeneous} of degree $d$ if every monomial appearing in it has total degree $d$. For a graded domain $R$, the graded analogue of a notion typically involves only the homogeneous elements and homogeneous ideals of $R$. The graded analogue of the quotient field $K$ is the \textit{homogeneous quotient field} $R_H = \{r/h : r \in R,\ h \in H\}$, where $H$ is the set of nonzero homogeneous elements of $R$. In particular, graded overrings of $R$ lie between $R$ and $R_H$ and have the inherited grading. For example, $R$ is a \textit{gr-valuation domain} if $x/y \in R$ or $y/x \in R$ for all $x, y \in H$ (see \cite{ACZ18} and \cite{AAC17}), and $R$ is a \textit{gr-Pr\"ufer domain} if every nonzero finitely generated homogeneous ideal of $R$ is invertible (see \cite{ACZ18}).

In \cite{Mer24}, Merchant introduces \textit{graded $i$-extensions} and \textit{graded $i$-domains}, and Hensler and Klawa \cite{HK25} further study these notions and their connection to gr-Pr\"ufer domains. In particular, they prove a graded analogue of \cite[Proposition~2.14]{Pap76}; that is, they show that a graded domain $R$ is a graded $i$-domain if and only if $R \subseteq \overline{R}$ is a graded $i$-extension and $\overline{R}$ is a gr-Pr\"ufer domain \cite[Proposition~4.1]{HK25}. We use this result repeatedly in Section~3.

Every graded domain that is an $i$-domain is a graded $i$-domain, since graded overrings are overrings and homogeneous prime ideals are prime ideals. However, the converse is not true, as Example~\ref{ex:gr-not-i} shows. Therefore, Papick's results do not apply directly to graded $i$-domains, which illustrates why these graded analogues are not trivial. On the other hand, if $R$ is trivially graded so that $R = R_0$, then every nonzero element of $R$ is homogeneous and $R_H = K$, so graded overrings, gr-valuation domains, and graded $i$-domains are simply overrings, valuation domains, and $i$-domains.

In this paper, we prove graded analogues of \cite[Propositions~2.18 and~2.20]{Pap76}. Section 2 states definitions and known results that we use for the proofs. In Section 3, we state the graded analogue of \cite[Corollary~2.15]{Pap76} from \cite{FM} as Proposition~\ref{cor:gr-local} and illustrate it in Example~\ref{ex:t2t3}. We then prove a graded version of \cite[Theorem~19.15 (6) $\Rightarrow$ (1)]{Gil72} in Lemma~\ref{lem:gr-gilmer}, which we use to prove the analogue of \cite[Proposition~2.18]{Pap76} in Proposition~\ref{prop:2.18}. Finally, we prove the analogue of \cite[Proposition~2.20]{Pap76} in Proposition~\ref{prop:2.20}.

\section{Preliminary Material}

In this section, we specify notation, definitions, and known results that we use. All rings are commutative with identity, and we refer to \cite{AM16} for standard background. 

A \textit{ring extension} $R \subseteq S$ is a ring $S$ together with a subring $R$. An element $s \in S$ is \textit{integral} over $R$ if it is a root of a monic polynomial in $R[X]$, and $S$ is \textit{integral} over $R$ if each of its elements is. If $R$ is an integral domain with quotient field $K$, an \textit{overring} of $R$ is a ring $S$ with $R \subseteq S \subseteq K$; for example, $\mathbb{Z}[1/2]$ and $\{a/b : a, b \in \mathbb{Z},\ b \text{ odd}\}$ are overrings of $\mathbb{Z}$. The \textit{integral closure} $\overline{R}$ of $R$ is the set of elements of $K$ that are integral over $R$. Furthermore, $\overline{R}$ is an overring of $R$, and $R$ is \textit{integrally closed} if $\overline{R} = R$. An integral domain $R$ with quotient field $K$ is \textit{local} if it has a unique maximal ideal, and a \textit{valuation domain} if $x \in R$ or $x^{-1} \in R$ for each nonzero $x \in K$. A \textit{valuation overring} of $R$ is an overring of $R$ that is a valuation domain.

For a ring $S$, let $\Spec(S)$ denote the set of prime ideals of $S$. For $a_1, \ldots, a_n \in S$, let $(a_1, \ldots, a_n)$ denote the ideal of $S$ generated by $a_1, \ldots, a_n$, so $(0) = \{0\}$ is the zero ideal. If $R \subseteq S$ is a ring extension, $P$ is a prime ideal of $R$, and $Q$ is a prime ideal of $S$ with $Q \cap R = P$, then we say that $Q$ \textit{lies over} $P$. As in \cite[p.~28]{Kap74}, $R \subseteq S$ satisfies \textit{lying over} (LO) if every prime of $R$ has a prime of $S$ lying over it. Similarly, $R \subseteq S$ satisfies \textit{incomparability} (INC) if distinct primes of $S$ lying over the same prime of $R$ are incomparable \cite[p.~28]{Kap74}. If $S$ is integral over $R$, then both hold \cite[Theorem~44]{Kap74}.

The definitions below for graded domains are standard (for example, see \cite[Section~2]{ACZ18} and \cite[Section~1.1]{SS26}). Let $\Gamma$ denote a torsionless grading monoid, which is a torsion-free commutative cancellative monoid. An integral domain $R$ is a \textit{$\Gamma$-graded domain} if
\[
R = \bigoplus_{\alpha \in \Gamma} R_\alpha,
\]
where each $R_\alpha$ is an additive subgroup of $R$, and $R_\beta R_\gamma \subseteq R_{\beta+\gamma}$ for all $\beta,\gamma \in \Gamma$. A nonzero $r \in R_\alpha$ is \textit{homogeneous} of degree $\alpha$. Each $r \in R$ can be written uniquely as $r = \sum_{\alpha} r_\alpha$, where $r_\alpha \in R_\alpha$ and only finitely many $r_\alpha$ are nonzero. We call the nonzero $r_\alpha$ the \textit{homogeneous components} of $r$. An ideal of $R$ is \textit{homogeneous} if it is generated by homogeneous elements.

\begin{exa}\label{ex:components}
Let $k$ be a field, and grade $R = k[X,Y]$ by total degree. That is, for each $n \geq 0$, $R_n$ consists of $0$ and the polynomials whose monomials all have total degree $n$, where the total degree of $cX^aY^b$ with $c \in k \setminus \{0\}$ is $a + b$. Then $X^2 + 1$ is not homogeneous, and its homogeneous components are $X^2$ and $1$. On the other hand, $X^2 + XY$ is homogeneous of degree $2$, since $X^2$ and $XY$ both have degree $2$, so its only homogeneous component is $X^2 + XY$ itself.
\end{exa}

Let $H = H(R)$ be the set of nonzero homogeneous elements of $R$; that is, $H = \bigcup_{\alpha \in \Gamma} (R_\alpha \setminus \{0\})$. Then $R_H = \{r/h : r \in R,\ h \in H\}$, the homogeneous quotient field of $R$, is a $\langle \Gamma \rangle = \{a - b : a, b \in \Gamma\}$-graded domain where
\[
(R_H)_\alpha = \left\{ \frac{f}{g} \mid f \in R_\beta,\ 0 \neq g \in R_\gamma, \text{ and } \beta - \gamma = \alpha \right\}
\]
for $\alpha \in \langle \Gamma \rangle$. A \textit{graded overring} of $R$ is an overring $S$ with $S \subseteq R_H$ and $S = \bigoplus_{\alpha \in \langle \Gamma \rangle} \bigl(S \cap (R_H)_\alpha\bigr)$; such an $S$ is itself a $\langle \Gamma \rangle$-graded domain. By \cite[Lemma 1.6]{Cha17}, the integral closure $\overline{R}$ is a graded overring of $R$; in particular, $\overline{R} \subseteq R_H$. We use this fact throughout this paper. For example, it makes $R \subseteq \overline{R}$ an extension of graded domains.

Recall that an extension $R \subseteq S$ is an \textit{injective extension} (or \textit{$i$-extension}) if the map $\Spec(S) \to \Spec(R)$ is injective, and $R$ is an \textit{$i$-domain} if $R \subseteq S$ is an $i$-extension for every overring $S$ of $R$. For a graded domain $T$, let $\grSpec(T)$ denote the set of homogeneous prime ideals of $T$. If $R \subseteq S$ is an extension of graded domains (that is, $R$ is a graded subring of $S$), then $Q \cap R$ is a homogeneous prime ideal of $R$ for each homogeneous prime ideal $Q$ of $S$. We call $R \subseteq S$ a \textit{graded $i$-extension} if the map $\grSpec(S) \to \grSpec(R)$ given by $Q \mapsto Q \cap R$ is injective. A graded domain $R$ is a \textit{graded $i$-domain} if $R \subseteq S$ is a graded $i$-extension for every graded overring $S$ of $R$ \cite{Mer24}. 

\begin{exa}\label{ex:grext}
Let $k$ be a field, and grade $k[X]$ and $k[X,Y]$ by total degree. Then $k[X] \subseteq k[X,Y]$ is an extension of graded domains. The homogeneous prime ideals $(Y)$ and $(0)$ of $k[X,Y]$ satisfy $(Y) \cap k[X] = (0) \cap k[X] = (0)$, so $k[X] \subseteq k[X,Y]$ is not a graded $i$-extension.
\end{exa}

Recall that $R = \bigoplus_{\alpha \in \Gamma} R_\alpha$ is a \textit{graded valuation domain} (gr-valuation domain) if either $\frac{x}{y} \in R$ or $\frac{y}{x} \in R$ for every $x,y \in H$. A \textit{graded valuation overring} of $R$ is a graded overring of $R$ that is a gr-valuation domain. By \cite[Theorem 2.2]{AAC17}, $R$ is a gr-valuation domain if and only if its homogeneous ideals are totally ordered under inclusion. Recall that $R$ is \textit{graded Pr\"ufer} (gr-Pr\"ufer) if every nonzero finitely generated homogeneous ideal of $R$ is invertible \cite{ACZ18}. Here, a nonzero ideal $I$ of $R$ is \textit{invertible} if $II^{-1} = R$, where $I^{-1} = \{x \in K : xI \subseteq R\}$ and $K$ is the quotient field of $R$.

\begin{exa}\label{ex:kX}
Let $k$ be a field, and grade $k[X]$ by degree. Then $k[X]$ is a gr-valuation domain, since its homogeneous elements have the form $aX^n$ with $a \in k \setminus \{0\}$. However, $k[X]$ is not a valuation domain, since neither $X/(X+1)$ nor $(X+1)/X$ lies in $k[X]$.
\end{exa}

The next example shows that a graded $i$-domain need not be an $i$-domain. It illustrates that Papick's results do not apply directly to graded $i$-domains.

\begin{exa}\label{ex:gr-not-i}
Let $k$ be a field, and grade $R = k[X, X^{-1}, Y]$ by $\Gamma = \mathbb{Z} \times \mathbb{N}_0$, where $\mathbb{N}_0 = \{0, 1, 2, \ldots\}$, so that each monomial $X^mY^n$ is homogeneous of degree $(m,n)$.  Then $R_H = R[Y^{-1}] = k[X, X^{-1}, Y, Y^{-1}]$, since $X$ is a unit of $R$ and the nonzero homogeneous elements of $R$ have the form $aX^mY^n$ with $a \in k \setminus \{0\}$ and $n \geq 0$. Moreover, every nonzero homogeneous element $aX^mY^n$ of $R_H$ is a unit. If $S \neq R$ is a graded overring of $R$ (where $R$ itself is a graded overring), then $S$ contains a homogeneous element $aX^mY^n \notin R$. Then $n < 0$, since $aX^mY^n \in R$ when $n \geq 0$. So, $a^{-1}X^{-m}Y^{-n-1} \in R$ since $-n-1 \geq 0$, and hence $Y^{-1} = (a^{-1}X^{-m}Y^{-n-1})(aX^mY^n) \in S$ and $S = R_H$. Since $(0)$ is the only homogeneous prime ideal of $R_H$, $R \subseteq R_H$ is a graded $i$-extension, and hence $R$ is a graded $i$-domain. 

However, $R$ is not an $i$-domain. To see this, let $u = Y/(X-1)$ and $T = R[u]$, which equals $k[X, X^{-1}, u]$ since $Y = (X-1)u$. Then $T$ is a polynomial ring in $u$ over $k[X, X^{-1}]$, since any polynomial relation for $u$ over $k[X, X^{-1}]$ would give a polynomial relation for $Y = (X-1)u$. Thus, for each $c \in k$, there is a ring homomorphism $\varphi_c \colon T \to k$ that fixes $k$ and sends $X$ to $1$ and $u$ to $c$, and its kernel $Q_c$ is a prime ideal of $T$. Since $\varphi_c(X) = 1$ and $\varphi_c(Y) = \varphi_c(X-1)\varphi_c(u) = 0$, the restriction of $\varphi_c$ to $R$ is $f \mapsto f(1,0)$, which does not depend on $c$. Hence, taking $c = 0$ and $c = 1$, we have $Q_0 \cap R = Q_1 \cap R$, but $u \in Q_0 \setminus Q_1$ since $\varphi_0(u) = 0$ and $\varphi_1(u) = 1$, so $R \subseteq T$ is not an $i$-extension.
\end{exa}

For an ideal $I$ of a graded domain $R$, let $I^*$ denote the ideal generated by the homogeneous elements of $I$; if $I$ is prime, then $I^*$ is a homogeneous prime ideal. A \textit{maximal homogeneous ideal} of $R$ is a proper homogeneous ideal that is maximal among proper homogeneous ideals. Every maximal homogeneous ideal of $R$ is prime \cite[Lemma~2.3]{ACZ18}, and every proper homogeneous ideal of $R$ is contained in at least one maximal homogeneous ideal \cite[Proposition~2.4]{ACZ18}. A graded domain $R$ is \textit{graded local} if it has a unique maximal homogeneous ideal. For $P \in \grSpec(R)$, $R_{H(R)\setminus P}$ denotes the localization of $R$ at $H(R)\setminus P$, which is a graded overring of $R$. 

An extension $R \subseteq S$ of graded domains satisfies \textit{graded lying over} if every homogeneous prime ideal of $R$ has a homogeneous prime ideal of $S$ lying over it, that is, if the map $\grSpec(S) \to \grSpec(R)$ is surjective. If $S$ is integral over $R$, then $R \subseteq S$ satisfies graded lying over. For $P \in \grSpec(R)$, lying over gives a prime ideal $Q$ of $S$ with $Q \cap R = P$ \cite[Theorem~44]{Kap74}. This $Q$ need not be homogeneous, but $Q^*$ is a homogeneous prime ideal of $S$ that still lies over $P$ (see the proof of \cite[Lemma~4.5]{SS26}).

\section{Graded Injective Domains}

We first state the graded analogue of \cite[Corollary~2.15]{Pap76}, which is established in \cite{FM}. We use this result in Example~\ref{ex:t2t3} and in Proposition~\ref{prop:2.20}.

\begin{prop}[{\cite{FM}}]\label{cor:gr-local}
A graded domain $R$ is a graded local graded $i$-domain if and only if $\overline{R}$ is a gr-valuation domain.
\end{prop}

\begin{exa}\label{ex:t2t3}
Let $k$ be a field, and grade $R = k[t^2, t^3]$ by degree in $t$. Since $t = t^3/t^2$ is a root of $X^2 - t^2 \in R[X]$, we have $t \in \overline{R}$, so the graded overring $\overline{R}$ of $R$ (being a ring containing $k$ and $t$) contains $k[t]$. The nonzero homogeneous elements of $R_H = k[t, t^{-1}]$ have the form $at^n$ with $a \in k \setminus \{0\}$ and $n \in \mathbb{Z}$, and for any $at^m$ and $bt^n$ of this form with $m \geq n$, we have $at^m/(bt^n) = (a/b)t^{m-n} \in k[t] \subseteq \overline{R}$. Hence, $\overline{R}$ is a gr-valuation domain, and $R$ is a graded local graded $i$-domain by Proposition~\ref{cor:gr-local}, although $R$ is neither integrally closed (as $t \notin R$) nor a gr-valuation domain (as neither $t^3/t^2$ nor $t^2/t^3$ lies in $R$). Note that \cite[Corollary~2.15]{Pap76} cannot give this conclusion, since $R$ is not local. This is because the nonunit $t^2 - 1$ of $R$ does not lie in the maximal ideal $(t^2, t^3)$.
\end{exa}

We now prove a graded analogue of the implication (6) $\Rightarrow$ (1) of \cite[Theorem 19.15]{Gil72}, which we need for Proposition~\ref{prop:2.18}.

\begin{lemma}\label{lem:gr-gilmer}
Let $D$ be an integrally closed graded domain, and let $P\in\grSpec(D)$. If each graded valuation overring of $D$ has at most one homogeneous prime ideal lying over $P$, then $D_P$ is a valuation domain.
\end{lemma}

\begin{proof} Let $H = H(D)$ be the set of nonzero homogeneous elements of $D$. We follow the methods in the proof of $({\sim}1)\Rightarrow({\sim}3)\Rightarrow({\sim}6)$ in \cite[Theorem 19.15]{Gil72}. Let $P$ be an arbitrary homogeneous prime ideal of $D$. We prove the contrapositive for $D_{H\setminus P}$: if $D_{H\setminus P}$ is not a gr-valuation domain, then some graded valuation overring of $D$ has two distinct homogeneous prime ideals lying over $P$. This suffices because if $D_{H\setminus P}$ is a gr-valuation domain, then $D_P$ is a valuation domain. To see why, note that a homogeneous element $a/s$ of $D_{H\setminus P}$, with $a \in H$ and $s \in H \setminus P$, lies in $PD_{H\setminus P}$ if $a \in P$ and is a unit otherwise. Hence, $PD_{H\setminus P}$ is the unique maximal homogeneous ideal of $D_{H\setminus P}$. Thus, if $D_{H\setminus P}$ is a gr-valuation domain, then $D_P=(D_{H\setminus P})_{PD_{H\setminus P}}$ is a valuation domain by \cite[Theorem 2.3(3),(5)]{AAC17}.

So, suppose $D_{H\setminus P}$ is not a gr-valuation domain. Then there are nonzero homogeneous $x,y\in D_{H\setminus P}$ with $x/y, y/x\notin D_{H\setminus P}$, so $t=x/y$ is a homogeneous element of $D_H$ such that neither $t$ nor $t^{-1}$ is in $D_{H\setminus P}$. Then neither $t$ nor $t^{-1}$ is in $D_P$. To see why, suppose $t\in D_P$, say $t=d/s$ with $d\in D$ and $s\in D\setminus P$, and choose a homogeneous component $s'$ of $s$ that is not in $P$. Since $t$ is homogeneous, $ts'$ is a homogeneous component of $ts=d$, so $ts'\in D$ and $t=ts'/s'\in D_{H\setminus P}$, a contradiction. The same argument applies to $t^{-1}$. Hence, by \cite[Lemma 19.14]{Gil72}, $t$ is not a root of any polynomial in $D[X]\setminus P[X]$. That is, as in the proof of $({\sim}1)\Rightarrow({\sim}3)$ in \cite[Theorem 19.15]{Gil72}, the kernel $Q$ of the canonical homomorphism $f(X)\mapsto f(t)$ of $D[X]$ onto $D[t]$ is contained in $P[X]$. Moreover, $Q$ is prime since $D[X]/Q \cong D[t]$ is a domain, and $Q \cap D = (0)$ since the homomorphism $f(X) \mapsto f(t)$ is the identity on $D$.

As in the proof of $({\sim}3)\Rightarrow({\sim}6)$ in \cite[Theorem 19.15]{Gil72}, the ideal $P_0=P[X]/Q=PD[t]$ of $D[t]\cong D[X]/Q$ satisfies
\[
D[t]/P_0\cong D[X]/P[X]\cong(D/P)[X].
\]
Since $(D/P)[X]$ is a domain, $P_0$ is prime. Let $M$ be the ideal of $D[t]$ generated by $P$ and $t$. Let $\psi \colon D[t] \to D/P$ be the homomorphism obtained by following the quotient map $D[t] \to D[t]/P_0 \cong (D/P)[X]$ with the map $f(X) \mapsto f(0)$. Since the isomorphism $D[t]/P_0 \cong (D/P)[X]$ sends each $d \in D$ to $d + P$ and $t$ to $X$, $\psi$ sends each element $d_0 + d_1t + \cdots + d_nt^n$ of $D[t]$, where each $d_i \in D$, to $d_0 + P$. In particular, $\psi$ is surjective and sends $P$ and $t$ to $0$, so $M \subseteq \ker(\psi)$. Conversely, if $\psi(d_0 + d_1t + \cdots + d_nt^n) = 0$, then $d_0 \in P$, so $d_0 + d_1t + \cdots + d_nt^n = d_0 + t(d_1 + d_2t + \cdots + d_nt^{n-1}) \in M$. Hence, $\ker(\psi) = M$, and $D[t]/M \cong D/P$ by the first isomorphism theorem. Since $D/P$ is a domain, $M$ is prime, and $M\neq P_0$ since $t\in M\setminus P_0$. Moreover, $P_0$ and $M$ both lie over $P$, which we use at the end of the proof. Since $P \subseteq P_0 \subseteq M$, we have $P \subseteq P_0 \cap D \subseteq M \cap D$, so it suffices to show that $M \cap D \subseteq P$. If $d \in M \cap D$, then $d + P = \psi(d) = 0$ since $\ker(\psi) = M$, so $d \in P$. Hence, $P_0 \cap D = M \cap D = P$.

Since $t$ is homogeneous, $D[t]$ is a graded overring of $D$ (see the remark preceding \cite[Proposition~3.3]{SS26}), and $P_0$ and $M$ are homogeneous, being generated by homogeneous elements. By \cite[Lemma 4.4]{CO18}, which plays the role of \cite[Corollary 19.7]{Gil72} in Gilmer's proof, there is a graded valuation overring $V$ of $D[t]$ containing homogeneous prime ideals $P_1$ and $M_1$ with $P_1\cap D[t]=P_0$ and $M_1\cap D[t]=M$. Since every nonzero homogeneous element of $D_H$ is a unit, $D$ and $D[t]$ have the same homogeneous quotient field $D_H$, so $V$ is also a graded valuation overring of $D$. Finally, $P_1\neq M_1$ since $P_0\neq M$, and they both lie over $P$ since $P_0$ and $M$ do as well. Thus, $V$ is a graded valuation overring of $D$ with two distinct homogeneous prime ideals lying over $P$.
\end{proof}

We now prove the graded analogue of \cite[Proposition 2.18]{Pap76}.

\begin{prop}\label{prop:2.18}
Let $R$ be a graded domain. Then $R$ is a graded $i$-domain if and only if $R \subseteq V$ is a graded $i$-extension for each graded valuation overring $V$ of $R$ and $R \subseteq \overline{R}$ is a graded $i$-extension.
\end{prop}

\begin{proof}

($\implies$) Suppose $R$ is a graded $i$-domain, and let $V$ be a graded valuation overring of $R$. By definition, since $V$ and $\overline{R}$ are graded overrings of $R$, $R \subseteq V$ and $R \subseteq \overline{R}$ are graded $i$-extensions.

($\impliedby$) Suppose $R \subseteq V$ is a graded $i$-extension for each graded valuation overring $V$ of $R$, and $R \subseteq \overline{R}$ is a graded $i$-extension. 

By \cite[Theorem 2.10]{Joh83}, $\overline{R}$ is the intersection of the graded valuation overrings of $R$; in particular, $\overline{R} \subseteq V$ for every graded valuation overring $V$ of $R$. Since $R \subseteq \overline{R} \subseteq R_H$, every graded valuation overring of $\overline{R}$ is a graded valuation overring of $R$. Let $V$ be an arbitrary graded valuation overring of $\overline{R}$. If $Q_1, Q_2 \in \grSpec(V)$ satisfy $Q_1 \cap \overline{R} = Q_2 \cap \overline{R}$, then $Q_1 \cap R = Q_2 \cap R$, so $Q_1 = Q_2$ because $R \subseteq V$ is a graded $i$-extension. Hence, $\overline{R} \subseteq V$ is a graded $i$-extension. 

In particular, for each $P \in \grSpec(\overline{R})$, at most one homogeneous prime ideal of each graded valuation overring of $\overline{R}$ lies over $P$. Thus, since $\overline{R}$ is integrally closed, $\overline{R}_P$ is a valuation domain for each such $P$ by Lemma~\ref{lem:gr-gilmer}, and hence $\overline{R}$ is gr-Pr\"ufer by \cite[Theorem 3.1((5) $\Rightarrow$ (1))]{ACZ18}. Since $R \subseteq \overline{R}$ is a graded $i$-extension and $\overline{R}$ is gr-Pr\"ufer, $R$ is a graded $i$-domain by \cite[Proposition~4.1]{HK25}.
\end{proof}

We now prove the graded analogue of \cite[Proposition 2.20]{Pap76}.

\begin{prop}\label{prop:2.20}
Let $R \subseteq T$ be an extension of graded domains such that $T$ is a graded $i$-domain integral over $R$. If $T$ is graded local or $R$ is integrally closed, then $R$ is a graded $i$-domain.
\end{prop}

\begin{proof} First, assume that $T$ is graded local, and let $S$ be a graded overring of $R$. By Proposition~\ref{cor:gr-local}, $\overline{T}$ is a gr-valuation domain, so its homogeneous prime ideals are totally ordered \cite[Theorem 2.2]{AAC17}. Since $T \subseteq \overline{T}$ satisfies graded lying over, every homogeneous prime ideal of $T$ is the contraction of a homogeneous prime ideal of $\overline{T}$. Since contractions of comparable ideals are comparable, any two homogeneous prime ideals of $T$ are comparable.

Moreover, $R \subseteq T$ is a graded $i$-extension. To see this, let $Q_1, Q_2 \in \grSpec(T)$ with $Q_1 \cap R = Q_2 \cap R$. By the above argument, $Q_1$ and $Q_2$ are comparable. Since $R \subseteq T$ is integral, it satisfies INC \cite[Theorem 44]{Kap74}, so distinct prime ideals of $T$ lying over the same prime of $R$ are incomparable. Hence, $Q_1 = Q_2$. As in the proof of \cite[Proposition 2.20]{Pap76}, consider the following diagram:
\[
\begin{tikzcd}
T \arrow[r, hook] & TS \\
R \arrow[u, hook] \arrow[r, hook] & S \arrow[u, hook]
\end{tikzcd}
\]
Here, $TS$ denotes the subring of the quotient field of $T$ generated by $T$ and $S$, which consists of the finite sums of products $ts$ with $t \in T$ and $s \in S$. It is a graded overring of $T$ with $TS = \bigoplus_{\gamma} (TS)_\gamma$, where $(TS)_\gamma = \sum_{\alpha+\beta=\gamma} T_\alpha S_\beta$ (see the proof of \cite[Lemma~4.5]{SS26}). Then $T \subseteq TS$ is a graded $i$-extension since $T$ is a graded $i$-domain. Since $R \subseteq T \subseteq TS$ and $R \subseteq T$ is also a graded $i$-extension, it follows that $R \subseteq TS$ is a graded $i$-extension. However, $TS$ is generated as a ring by $S$ and $T$, and every element of $T$ is integral over $R \subseteq S$, so $TS$ is integral over $S$ \cite[Corollary~5.3]{AM16}. Hence, $S \subseteq TS$ satisfies graded lying over. Thus, if $P_1, P_2 \in \grSpec(S)$ satisfy $P_1 \cap R = P_2 \cap R$, there are $P_1', P_2' \in \grSpec(TS)$ with $P_1' \cap S = P_1$ and $P_2' \cap S = P_2$. Then $P_1' \cap R = P_2' \cap R$, so $P_1' = P_2'$, and hence $P_1 = P_2$. Therefore, $R \subseteq S$ is a graded $i$-extension, and $R$ is a graded $i$-domain.

Now suppose that $R$ is integrally closed. By \cite[Proposition 4.1]{HK25}, $\overline{T}$ is gr-Pr\"ufer. Let $P \in \grSpec(R)$. Since $\overline{T}$ is integral over $R$ \cite[Theorem 40]{Kap74}, there is $N \in \grSpec(\overline{T})$ with $N \cap R = P$. Since $\overline{T}_N$ is a valuation domain \cite[Theorem 3.1((1)$\Rightarrow$(5))]{ACZ18}, $R_P$ is a valuation domain by \cite[Theorem 22.4]{Gil72}. Thus, $R$ is gr-Pr\"ufer \cite[Theorem 3.1((5)$\Rightarrow$(1))]{ACZ18}. Since $R = \overline{R}$, $R$ is a graded $i$-domain by \cite[Proposition 4.1]{HK25}. 
\end{proof}

\section*{Acknowledgements}

I am very grateful for the guidance and feedback provided by Dr. Hannah Klawa in the writing of this paper.

\bibliographystyle{amsalpha}
\bibliography{bibliography}

\end{document}